\documentclass{amsart}

\newtheorem{theorem}{Theorem}[section]
\newtheorem{lemma}[theorem]{Lemma}
\newtheorem{corollary}[theorem]{Corollary}
\newtheorem{proposition}[theorem]{Proposition}

\theoremstyle{definition}

\newtheorem{example}[theorem]{Example}

\theoremstyle{remark}
\newtheorem{remark}[theorem]{Remark}

\numberwithin{equation}{section}

\newcommand{\dis}{\displaystyle}

\begin{document}

\title[Set-theoretic complete
intersection monomial curves] {PRODUCING SET-THEORETIC COMPLETE
INTERSECTION MONOMIAL CURVES IN $\mathbb{P}^n$}
\author{MESUT \c{S}AH\.{I}N}

\address{Department of Mathematics, Bilkent University,
06800 Ankara, TURKEY} \email{sahin@fen.bilkent.edu.tr}

\curraddr{Department of Mathematics, At\i l\i m University, 06836
Ankara, TURKEY} \email{mesut@atilim.edu.tr}


\date{\today}


\maketitle

\begin{abstract}
In this paper we describe an algorithm for  producing infinitely
many examples of set-theoretic complete intersection monomial
curves in $\mathbb{P}^{n+1}$, starting with a single set-theoretic
complete intersection  monomial curve in $\mathbb{P}^{n}$.
Moreover we investigate the numerical criteria to decide when
these monomial curves can or cannot be obtained via semigroup
gluing.
\end{abstract}

\section{Introduction}

It is well known that a variety in an $n$-space can be written as
the intersection of $n$ hypersurfaces set theoretically, see
\cite{eis}. It is then natural to ask whether this number is
minimal. A curve in  $n$-space which is the intersection of $n-1$
hypersurfaces is called a set-theoretic complete intersection,
s.t.c.i. for short.  If moreover its defining ideal is generated
by $n-1$ polynomials, then it is called an ideal theoretic
complete intersection, abbreviated i.t.c.i.. Determining
set-theoretic or ideal-theoretic complete intersection curves is a
classical and longstanding problem in algebraic geometry. An
associated problem is to give explicitly the equations of the
hypersurfaces involved. When the characteristic of the field K is
positive, it is known that all monomial curves are s.t.c.i. in
$\mathbb{P}^{n}$, see \cite{moh}. However the question is still
open in characteristic zero case despite the tremendous progress
in this direction, see for example \cite{eto2,kat,thoma4} and the
references there for some recent activity.

The purpose of the present paper is to describe a  method to
produce \textit{infinitely many} s.t.c.i. monomial curves starting
from one single s.t.c.i. monomial curve, see section 4. Our
approach has the side novelty of describing explicitly the
equations of hypersurfaces on which these new monomial curves lie
as s.t.c.i.. On the other hand, semigroup gluing being one of the
most popular techniques of recent research, we develop numerical
criteria to determine when these new curves can or cannot be
obtained via gluing, see section 3. In the last section we discuss
several consequences and variations of these results.

\section{Preliminaries}

Throughout the paper, $K$ will be assumed to be an algebraically
closed field of characteristic zero. By an \textit{affine monomial
curve} $C(m_{1},\dots,m_{n})$, for some positive integers
$m_{1}<\cdots<m_{n}$ with $gcd(m_{1},\dots,m_{n})=1$, we mean a
curve with generic zero $(v^{m_{1}},\dots ,v^{m_{n}})$ in the
affine n-space $\mathbb{A}^{n}$, over $K$. By a \textit{projective
monomial curve} $\overline{C}(m_{1},\dots,m_{n})$ we mean a curve
with generic zero
\[(u^{m_{n}},u^{m_{n}-m_{1}}v^{m_{1}},\dots
,u^{m_{n}-m_{n-1}}v^{m_{n-1}},v^{m_{n}})
\]
in the projective n-space $\mathbb{P}^{n}$, over $K$. Note that
$\overline{C}(m_{1},\dots,m_{n})$ is the projective closure of
$C(m_{1},\dots,m_{n})$.

Whenever we write $\overline{C} \subset \mathbb{P}^n$ to simplify
the notation, we always mean a monomial curve
$\overline{C}(m_{1},\dots,m_{n})$ for some fixed positive integers
$m_1<\dots<m_n$ with $gcd(m_1,\dots,m_n)=1$.

Let $m$ be a positive integer in the numerical semigroup generated
by $m_1,\dots,m_n$, i.e. $m=s_1 m_1 +\cdots+ s_n m_n$ where
$s_1,\dots ,s_n$ are some non-negative integers. Note that in
general there is no unique choice for $s_1,\dots,s_n$ to represent
$m$ in terms of $m_1,\dots,m_n$. We define the degree $\delta(m)$
of $m$ to be the minimum of all possible sums $s_1+\cdots+s_n$. If
$\ell$ is a positive integer with $gcd(\ell,m)=1$, then we say
that the monomial curve $\overline{C}(\ell m_1,\dots,\ell m_n, m)$
in $\mathbb{P}^{n+1}$ is an {\em extension} of $\overline{C}$. We
similarly define $C(\ell m_1,\dots,\ell m_n, m)$ to be an {\em
extension} of $C=C(m_{1},\dots,m_{n})$. We say that  an extension
is {\em nice} if
 $\delta(m) > \ell$ and {\em bad} otherwise,
adopting the terminology of \cite{pf}.

When the integers $m_1,\dots,m_n$ are fixed and understood in a
discussion, we will use $\overline{C}_{\ell,m}$ to denote the
extensions $\overline{C}(\ell m_1,\dots,\ell m_n, m)$ in
$\mathbb{P}^{n+1}$, and use $C_{\ell,m}$ to denote the extensions
$C(\ell m_1,\dots,\ell m_n, m)$ in $\mathbb{A}^{n+1}$.

\subsection{Extensions of Monomial Curves in $\mathbb{A}^{n}$}

Let $C=C(m_{1},\dots,m_{n})$ be a s.t.c.i. monomial curve in
$\mathbb{A}^{n}$. In this section, we show that all extensions of
$C$, in the sense defined above, are s.t.c.i. For this we first
define, for any ideal $I \subset K[x_1,\dots,x_{n+1}]$,
$\Gamma_{\ell}(I)$ to be the ideal which is generated by all
polynomials of the form $\Gamma_{\ell}({g})$, where
$\Gamma_{\ell}({g(x_1,\dots ,x_{n+1})})=g(x_1,\dots
,x_n,x_{n+1}^{\ell})$, for all $g\in I$. We use the following
trick of M. Morales:

\begin{lemma}[\mbox{\cite[Lemma ~3.2]{mor}}]\label{mor} Let $Y_{\ell}$ be the monomial curve $C(\ell m_1,\dots,\ell m_n,m_{n+1})$ in
$\mathbb{A}^{n+1}$. Then $I(Y_{\ell})=\Gamma_{\ell}(I(Y_1))$.
\end{lemma}

For any extension of $C$ of the form $C_{\ell,m}$, we obviously
have $I(C)\subset I(C_{\ell,m})$ and $I(C_{\ell,m})\cap
K[x_1,\dots,x_n]=I(C)$. The exact relation between the ideals of
$C$ and $C_{\ell,m}$ are given by the following lemma.

\begin{lemma}\label{afinlemma} Let $m=s_1 m_1 +\cdots+ s_n m_n$. For any positive integer $\ell$ with $gcd(\ell,m)=1$ we have
$I(C_{\ell,m})=I(C)+(G)$, where $G={x_{1}}^{s_{1}}\cdots
\,\,{x_{n}}^{s_{n}}-x_{n+1}^{\ell}$.
\end{lemma}

\begin{proof} \mbox{} \par
 \textbf{Case $\ell=1$:}  We  show that
$I(C_{1,m})=I(C)+({x_{1}}^{s_{1}}\cdots \,\,{x_{n}}^{s_{n}}-x_{n+1})$.\\
For any polynomial $f\in K[x_{1},\dots ,x_{n+1}]$, there are
polynomials $g\in K[x_{1},\dots ,x_{n}]$ and $h\in K[x_{1},\dots
,x_{n+1}]$ such that
\begin{eqnarray*}
f(x_{1},\dots ,x_{n+1})&=& f(x_{1},\dots
,x_{n},x_{n+1}-x_{1}^{s_{1}}\cdots x_{n}^{s_{n}}+
x_{1}^{s_{1}}\cdots x_{n}^{s_{n}}) \\
&=&g(x_{1},\dots ,x_{n})+(x_{1}^{s_{1}}\cdots
x_{n}^{s_{n}}-x_{n+1})h(x_{1},\dots ,x_{n+1}).
\end{eqnarray*}
 This identity  implies that
$f\in I(C_{1,m})$ if and only if $g\in I(C)$.

\textbf{Case $\ell>1$:} Applying Lemma \ref{mor} with
$Y_1=C_{1,m}$ we have
\begin{eqnarray*}
I(C_{\ell,m})&=&\Gamma_{\ell}(I(C_{1,m})), \; \; \text{by Lemma \ref{mor}} \\
&=& \Gamma_{\ell}(I(C)+(x_{1}^{s_{1}}\cdots
x_{n}^{s_{n}}-x_{n+1})) \; \; \text{by the first part of this lemma} \\
&=& I(C)+(G).
\end{eqnarray*}
\end{proof}

This lemma provides an alternate proof to the following theorem
which is a special case of \cite[Theorem~2]{thoma4}.

\begin{theorem}\label{afinstci}
If{~} $C \subset \mathbb{A}^{n}$ is a s.t.c.i. monomial curve,
then all extensions of the form $C_{\ell,m} \subset
\mathbb{A}^{n+1}$ are also s.t.c.i. monomial curves.
\end{theorem}

\begin{proof} Since
$I(C_{\ell,m})=I(C)+(G)$ by Lemma \ref{afinlemma}, it follows that
\begin{eqnarray*}
Z(I(C_{\ell,m}))&=& Z(I(C)+(G)) \\
 C_{\ell,m}&=&Z(I(C))\bigcap Z(G),
\end{eqnarray*}
where $Z(\cdot)$ denotes the zero set as usual. Hence
 $C_{\ell,m}$ is a s.t.c.i. if $C$ is.
\end{proof}

\section{Extensions that can not be obtained by gluing}

If $\overline{C}(m_1,\dots ,m_{n+1})$ is a monomial curve in
$\mathbb{P}^{n+1}$, then there is a corresponding semigroup
$\mathbb{N}T$, where
\[ T=\{(m_{n+1},0),(m_{n+1}-m_1,m_1),\dots ,(m_{n+1}-m_n,m_n),(0,m_{n+1})\} \subset \mathbb{N}^2. \]

Let $T=T_1\bigsqcup T_2$ be a decomposition of $T$ into two
disjoint proper subsets. Without loss of generality assume that
the cardinality of $T_1$ is less than or equal to the cardinality
of $T_2$. $\mathbb{N}T$ is called a \textit{gluing} of
$\mathbb{N}T_1$ and $\mathbb{N}T_2$ if there exists a nonzero
$\alpha \in \mathbb{N}T_1 \bigcap \mathbb{N}T_2$ such that
$\mathbb{Z}\alpha=\mathbb{Z}T_1 \bigcap \mathbb{Z}T_2$. Following
the literature we write $I(T)$ for the ideal of the toric variety
corresponding to the affine semigroup $\mathbb{N}T$. Note that if
$\mathbb{N}T$ is a gluing of $\mathbb{N}T_1$ and $\mathbb{N}T_2$
then we have $I(T)=I(T_1)+I(T_2)+(G_\alpha)$, where $G_\alpha$ is
the relation polynomial, see \cite{thoma4}.

We note that the condition $\mathbb{Z}\alpha=\mathbb{Z}T_1 \bigcap
\mathbb{Z}T_2$ is not fulfilled when $T_1$ is not a singleton.
Hence we formulate this observation to be the following
\begin{proposition}\label{gluingprop1} If $T_1$ is not a singleton then $\mathbb{N}T$ is not a gluing of $\mathbb{N}T_1$ and
$\mathbb{N}T_2$.
\end{proposition}
\begin{proof} If $T_1$ is not a singleton, then neither is $T_2$ by the
assumption on the cardinalities of these sets. Thus
$\mathbb{Z}T_1$ and $\mathbb{Z}T_2$ are submodules of
$\mathbb{Z}^2$ of rank two each. It is elementary to show that
their intersection has rank two. For instance, let $r$ and $t$ be
generators of $\mathbb{Z}T_1$, then the images of $r$ and $t$ have
finite order in the finite group $\mathbb{Z}^2/\mathbb{Z}T_2$,
meaning that $ar$ and $bt$ are in $\mathbb{Z}T_2$ for some
positive integers $a$ and $b$. Then the rank two
$\mathbb{Z}$-module generated by $ar$ and $bt$ is contained in the
intersection $\mathbb{Z}T_1\cap \mathbb{Z}T_2$ which must be of
rank two itself being a submodule of $\mathbb{Z}^2$.

 Hence  the intersection cannot be generated by a single element. Thus
$\mathbb{N}T$ is not a gluing of $\mathbb{N}T_1$ and
$\mathbb{N}T_2$.
\end{proof}

This proposition means that the only way to show that an extension
in $\mathbb{P}^{n+1}$ is a s.t.c.i. via gluing is to apply the
technique to a projective monomial curve in $\mathbb{P}^{n}$. Thus
we discuss the case where $T_1$ is a singleton. But if $T_1$ is
$\{(m_{n+1},0)\}$ or $\{(0,m_{n+1})\}$ then $\mathbb{N}T_1 \bigcap
\mathbb{N}T_2=\{(0,0)\}$. So it is sufficient to deal with the
case where $T_1$ is of the form $\{(m_{n+1}-m_{i},m_{i})\}$, for
some $i \in \{1,\dots ,n\}$.

From now on, $\Delta_i$ denotes the greatest common divisor of the
positive integers $m_1,\dots ,\widehat{m_i},\dots ,m_{n+1}$ ($m_i$
is omitted), for  $i=1,\dots,n$. Note that we have
$gcd(\Delta_i,m_i)=1$, for all $i=1,\dots,n$, since
$gcd(m_1,\dots,m_{n+1})=1$.

\begin{proposition}\label{gluingprop2} If \;$T_1=\{(m_{n+1}-m_{i_0},m_{i_0})\}$ for
some fixed $i_0 \in \{1,\dots ,{n}\}$, then $\mathbb{N}T$ is a
gluing of $\mathbb{N}T_1$ and $\mathbb{N}T_2$
if and only if there exist non-negative integers $d_j$, for $j=1,\dots,\widehat{i}_0,\dots,n+1$, satisfying the following two conditions: \\
{\rm (I)} $\dis \Delta_{i_0}m_{i_0}=\sum_{\substack{j=1 \\j \neq
i_0}}^{n+1} d_jm_j$, and {\rm (II)}  $\dis \Delta_{i_0} \geq
\sum_{\substack{j=1 \\j \neq i_0}}^{n+1} d_j$.
\end{proposition}

\begin{proof} Let $\alpha=\Delta_{i_0}(m_{n+1}-m_{i_0},m_{i_0})$.
We first show that $\mathbb{Z}T_1 \bigcap
\mathbb{Z}T_2=\mathbb{Z}\alpha$. Since $\Delta_{i_0}=gcd(m_1,\dots
,\widehat{m_{i_0}},\dots ,m_{n+1})$, there are $z_j \in
\mathbb{Z}$, for $j=1,\dots,\widehat{i}_0,\dots,n+1$, such that
$\Delta_{i_0}=\sum_{j \neq i_0} z_jm_j$. So,
$\Delta_{i_0}m_{i_0}=\sum_{j \neq {i_0}} m_{i_0}z_jm_j$ which
implies that
\[ \Delta_{i_0}(m_{n+1}-m_{i_0},m_{i_0})=\sum_{j \neq {i_0}} m_{i_0}z_j(m_{n+1}-m_j,m_j)+(\Delta_{i_0}-\sum_{j \neq
{i_0}} m_{i_0}z_j)(m_{n+1},0).\]
 Thus $\alpha=\Delta_{i_0}(m_{n+1}-m_{i_0},m_{i_0})\in
\mathbb{Z}T_1 \bigcap \mathbb{Z}T_2$ implying $\mathbb{Z}\alpha
\subseteq \mathbb{Z}T_1 \bigcap \mathbb{Z}T_2$.

For the converse inclusion, take $c(m_{n+1}-m_{i_0},m_{i_0})\in
\mathbb{Z}T_1 \bigcap \mathbb{Z}T_2$, for some $c\in \mathbb{Z}$.
Then, obviously we have $c(m_{n+1}-m_{i_0},m_{i_0}) \in
\mathbb{Z}T_2$ which implies that $cm_{i_0} \in
\mathbb{Z}(\{m_1,\dots ,\widehat{m_{i_0}},\dots
,m_{n+1}\})=\mathbb{Z}\Delta_{i_0}$. So, $\Delta_{i_0}$ divides
$cm_{i_0}$. If $\Delta_{i_0}>1$, then $\Delta_{i_0}$ divides $c$,
since it does not divide $m_{i_0}$ (remember that
$gcd(\Delta_{i_0},m_{i_0})=1$). If $\Delta_{i_0}=1$, obviously
$\Delta_{i_0}$ divides $c$. Thus, $c(m_{n+1}-m_{i_0},m_{i_0})$ is
a multiple of $\alpha$ and $\mathbb{Z}T_1 \bigcap \mathbb{Z}T_2
\subseteq \mathbb{Z}\alpha$.

Since $\mathbb{Z}T_1 \bigcap \mathbb{Z}T_2 = \mathbb{Z}\alpha$, it
will follow by definition that $\mathbb{N}T$ is a gluing of
$\mathbb{N}T_1$ and $\mathbb{N}T_2$ if and only if $\alpha \in
\mathbb{N}T_1 \bigcap \mathbb{N}T_2$. But, if $\alpha \in
\mathbb{N}T_1 \bigcap \mathbb{N}T_2$ then there exists
non-negative integers $d_j$ and $d$ for which we have
\begin{eqnarray*}\Delta_{i_0}(m_{n+1}-m_{i_0},m_{i_0})&=&\sum_{j \neq {i_0}}
d_j(m_{n+1}-m_j,m_j)+d(m_{n+1},0)\\
(\Delta_{i_0}m_{n+1}-\Delta_{i_0}m_{i_0},\Delta_{i_0}m_{i_0})
&=&([d+\sum_{j \neq {i_0}} d_j]m_{n+1}-\sum_{j \neq {i_0}}
d_jm_j,\sum_{j \neq {i_0}} d_jm_j).
\end{eqnarray*}
Thus, $\Delta_{i_0}m_{i_0}=\sum_{j \neq {i_0}} d_jm_j$ and
$d=\Delta_{i_0}-\sum_{j \neq {i_0}} d_j$. Since $d\geq0$, we see
that the conditions ${\rm (I)}$ and ${\rm (II)}$ hold. On the
other hand, if ${\rm (I)}$ and ${\rm (II)}$ hold then we observe
that $\alpha \in \mathbb{N}T_1 \bigcap \mathbb{N}T_2$, by the
equalities above. Thus, the condition $\alpha \in \mathbb{N}T_1
\bigcap \mathbb{N}T_2$ is equivalent to the existence of the
non-negative integers $d_j$ satisfying ${\rm (I)}$ and ${\rm
(II)}$.
\end{proof}

As a direct consequence of Proposition \ref{gluingprop2} we get
the following

\begin{corollary}\label{gluingcor1} If $\Delta_{i_0}=1$, for some fixed $i_0 \in
\{1,\dots ,{n}\}$, then $\mathbb{N}T$ cannot be obtained as a
gluing of $\mathbb{N}T_1$ and $\mathbb{N}T_2$, where
$T_1=\{(m_{n+1}-m_{i_0},m_{i_0})\}$ and $T_2=T - T_1$.

\end{corollary}

\begin{proof} We apply Proposition \ref{gluingprop2}.
If ${\rm (I)}$ does not hold, we are done. If it holds, then we
have two cases: either $\dis \sum_{\substack{j=1 \\j \neq
i_0}}^{n+1} d_j =1$ or $\dis \sum_{\substack{j=1 \\j \neq
i_0}}^{n+1} d_j
>1$. The first case forces $m_{i_0}=m_j$ for some $j\neq i_0$, from (I), but this contradicts the way we choose $m_i's$. The second case causes (II) to fail, as $\Delta_{i_0}=1$.
\end{proof}

\begin{example} If we consider the curve
$\overline{C}(2,3,4,8)\subset \mathbb{P}^4$ and take $i_0=2$, then
the conditions ${\rm (I)}$ and ${\rm (II)}$ of the above
proposition hold. Thus this curve can be obtained by gluing.

But if we consider the monomial curve
$\overline{C}(2,4,7,8)\subset \mathbb{P}^4$, then for every choice
of $i_0$, either $\Delta_{i_0}=1$, or else condition  ${\rm (II)}$
of the above proposition fails. Hence this curve cannot be
obtained by gluing.
\end{example}

\begin{corollary}\label{badextensions} Let $\overline{C}_{\ell,m}\subset \mathbb{P}^{n+1}$ be a bad
extension of $\overline{C}=\overline{C}(m_1,\dots ,m_n)$,  i.e.
$\ell \geq \delta(m)$. If $\overline{C}$ is a s.t.c.i. on the
hypersurfaces $f_1=\cdots=f_{n-1}=0$, then $\overline{C}_{\ell,m}$
can be shown to be a s.t.c.i. on the hypersurfaces
$f_1=\cdots=f_{n-1}=0$ and
$F=x_{n+1}^{\ell}-x_0^{\ell-\delta(m)}x_1^{s_1} \cdots
x_n^{s_n}=0$ by the technique of gluing, where
$m=s_1m_1+\cdots+s_nm_n$ and $s_1+\cdots+s_n=\delta(m)$.
\end{corollary}
\begin{proof}Since $m_1<\cdots<m_n$ and $m=s_1m_1+\cdots+s_nm_n \leq
\delta(m)m_n \leq \ell m_n$, it follows that $\ell m_n$ is the
biggest number among $\{\ell m_1,\dots ,\ell m_n,m\}$. The
extension $\overline{C}_{\ell,m}$ corresponds to the semigroup
$\mathbb{N}T$, where $T=T_1 \bigcup T_2$, $T_1=\{(\ell m_n-m,m)\}$
and $T_2=\{(\ell m_n,0),(\ell m_n-\ell m_1,\ell m_1),\dots ,(\ell
m_n-\ell m_{n-1},\ell m_{n-1}),(0,\ell m_n)\}$. Since $gcd(\ell
m_1,\dots ,\ell m_n)=\ell$, $\ell m=s_1(\ell m_1)+\cdots+s_n(\ell
m_n)$ and $\ell \geq \delta(m)$, $\mathbb{N}T$ is a gluing of
$\mathbb{N}T_1$ and $\mathbb{N}T_2$, by Proposition
\ref{gluingprop2}. Since $I(T)=I(T_1)+I(T_2)+(F)$, the claim
follows from \cite[Theorem~2]{thoma4}.
\end{proof}

\section{The Main Results}

Since \emph{bad} extensions are shown to be a s.t.c.i. by the
technique of gluing (see Corollary \ref{badextensions} above), we
study \textit{nice} extensions of monomial curves in this section.
By using the theory developed in the previous section one can
check which of these extensions can be obtained by the technique
of gluing semigroups.

Throughout this section we will assume that
\begin{itemize}
    \item $\overline{C}=\overline{C}(m_1,\dots,m_n) \subset \mathbb{P}^n$ is
a s.t.c.i. on $f_1=\cdots=f_{n-1}=0$
    \item $m=s_1m_1+\cdots+s_nm_n$ for some nonnegative integers
$s_1,\dots,s_n$ such that $s_1+\cdots+s_n=\delta(m)$
    \item $\ell$ is a positive integer with $gcd(\ell,m)=1$
    \item $\delta(m)>\ell$.
\end{itemize}

\begin{remark}\label{inclusion} Since $\overline{C}$ is s.t.c.i. on
$f_1=\cdots=f_{n-1}=0$, its affine part $C$ is s.t.c.i. on
$g_1=\cdots=g_{n-1}=0$, where $g_i(x_1,\dots ,x_n)=f_i(1,x_1,\dots
,x_n)$ is the dehomogenization of $f_i$,  $i=1,\dots ,n-1$. It
follows from Theorem \ref{afinstci} that $C_{\ell,m}$ is a
s.t.c.i. on the hypersurfaces $g_i=0$ and
$G={x_{1}}^{s_{1}}\cdots\,\,{x_{n}}^{s_{n}}-x_{n+1}^{\ell}=0$. So,
the ideal of the affine curve $C_{\ell,m}$ contains $g_i$'s and
$G$. Hence the ideal of the projective closure of $C_{\ell,m}$
must contain (at least) $f_i$'s and $F$, where $F$ is the
homogenization of $G$. Now, since $f_1,\dots , f_{n-1}, F \in
I(\overline{C}_{\ell,m})$, we always have $\overline{C}_{\ell,m}
\subseteq Z(f_1,\dots ,f_{n-1},F)$.
\end{remark}
\subsection{The case where $f_i$'s are
general, but $m$ is special} In this section we assume that $m$ is
a multiple of $m_n$, i.e. $m=s_nm_n$ where $s_n$ is a positive
integer. Note that $(s_1,\dots,s_{n-1})=(0,\dots,0)$ and
$\delta(m)=s_n$ in this case.

\begin{theorem}\label{projstci} Let $\overline{C}\subset \mathbb{P}^{n}$ be a s.t.c.i.
on the hypersurfaces $f_1=\cdots=f_{n-1}=0$, $gcd(\ell,s_nm_n)=1$
and $s_n>\ell$. Then, the nice extensions
$\overline{C}_{\ell,s_nm_n}$ in $\mathbb{P}^{n+1}$ are s.t.c.i. on
$f_1=\cdots=f_{n-1}=F=0$ where
$F=x_{n}^{s_n}-x_{0}^{s_n-\ell}x_{n+1}^{\ell}$.

\end{theorem}

\begin{proof} The fact that these nice extensions are s.t.c.i. can be seen easily by \cite[Theorem~3.4]{thoma2}
taking $b_1=m_1,\dots,b_{n-1}=m_{n-1}$, $d=m_n$ and
$k=(s_n-\ell)m_n$. In addition to this, we provide here the
equation of the binomial hypersurface $F=0$ on which these
extensions lie as s.t.c.i. monomial curves.

Since $\overline{C}_{\ell,s_nm_n} \subseteq Z(f_{1},\dots
,f_{n-1},F)$, we need to get the converse inclusion. Take a point
$P=(p_0,\dots ,p_n,p_{n+1}) \in Z(f_1,\dots ,f_{n-1},F)$. Then,
since $f_i \in K[x_0,\dots,x_n]$, we have $f_i(P)=f_i(p_0,\dots
,p_n)=0$, for all $i=1,\dots ,n-1$. Since $Z(f_1,\dots
,f_{n-1})=\overline{C}$ in $\mathbb{P}^n$ by assumption, the last
observation implies that
\[(p_0,\dots,p_n)=(u^{m_{n}},u^{m_{n}-m_{1}}v^{m_{1}},
\dots,u^{m_{n}-m_{n-1}}v^{m_{n-1}},v^{m_{n}}).
\]

If $p_0=0$ then $u=0$, yielding that
$(p_0,\dots,p_{n-1},p_n)=(0,\dots ,0,p_n)$. Since $s_n>\ell$, we
have also $p_n=0$, by $F(0,\dots
,0,p_n,p_{n+1})={p_{n}}^{s_n}-p_{0}^{s_n-\ell}p_{n+1}^{\ell}=0$.
So we observe that $(p_0,\dots ,p_n,p_{n+1})=(0,\dots ,0,1)$ which
is on the curve $\overline{C}_{\ell,s_nm_n}$. If $p_0=1$ then
$(1,p_1,\dots ,p_n,p_{n+1}) \in Z(g_1,\dots ,g_{n-1},G)$ by the
assumption, where $g_i$ and $G$ are polynomials defined in Remark
\ref{inclusion}. Since $C_{\ell,s_nm_n}$ is a s.t.c.i. on the
hypersurfaces $g_1=\cdots=g_{n-1}=0$ and $G=0$ it follows that
$(1,p_1,\dots ,p_n,p_{n+1}) \in C_{\ell,s_nm_n} \subset
\overline{C}_{\ell,s_nm_n}$.
\end{proof}

\subsection{The case where $f_i$'s are special and
$m$ is general} Assume now that $m$ is not a multiple of $m_n$,
i.e. $(s_1,\dots ,s_{n-1}) \neq (0,\dots ,0)$. Recall that we
choose $s_1,\dots,s_n$ in the representation of
$m=s_1m_1+\cdots+s_nm_n$ in such a way that $s_{1}+\dots +s_{n}$
is minimum, i.e. $s_{1}+\dots +s_{n}=\delta(m)$. First we prove a
lemma where no restriction on the $f_i$ is required.
\begin{lemma}\label{projlemma}
Let $\overline{C}\subset \mathbb{P}^{n}$ be a s.t.c.i.
 on $f_1=\cdots=f_{n-1}=0$ and $\delta(m)>\ell$. Then,
$Z(f_{1},\dots ,f_{n-1},F)= \overline{C}_{\ell,m}\cup L \subset
\mathbb{P}^{n+1}$, where
$F={x_{1}}^{s_{1}}\cdots\,\,{x_{n}}^{s_{n}}-x_{0}^{\delta(m)-
\ell}x_{n+1}^{\ell}$ and $L$ is the line $x_{0}=\cdots=x_{n-1}=0$.
\end{lemma}

\begin{proof} We first prove $\overline{C}_{\ell,m}\bigcup L \subseteq Z(f_1,\dots
,f_{n-1},F)$. By the light of Remark \ref{inclusion}, it is
sufficient to see that $L \subseteq Z(f_1,\dots ,f_{n-1},F)$. For
this, we take a point $P=(p_0,\dots ,p_{n+1})$ on the line $L$,
i.e., $P=(0,\dots ,0,p_n,p_{n+1})$. Since $(s_1,\dots ,s_{n-1})
\neq (0,\dots ,0)$ and $\delta(m)>\ell$, we see that $F(P)=0$.
Letting $v\in K$ be any $m_n$-th root of $p_n$, we get $(0,\dots
,0,p_n)=(0,\dots ,0,v^{m_n}) \in \overline{C}=Z(f_{1},\dots
,f_{n-1})$. Since the polynomials $f_i$ are in $K[x_0,\dots
,x_n]$, it follows that $f_i(P)=f_i(0,\dots ,0,p_n)=0$, for all
$i=1,\dots,n-1$. Thus $P \in Z(f_1,\dots ,f_{n-1},F)$.

For the converse inclusion, take $P=(p_0,\dots ,p_n,p_{n+1}) \in
Z(f_1,\dots ,f_{n-1},F)$. Then, for all $i=0,\dots ,n-1$, we get
$f_i(p_0,\dots ,p_n)=f_i(P)=0$ implying that
\[(p_0,\dots,p_n)=(u^{m_{n}},u^{m_{n}-m_{1}}v^{m_{1}},\dots,
u^{m_{n}-m_{n-1}}v^{m_{n-1}},v^{m_{n}}).
\]

If $p_0=0$ then $u=0$, yielding that $(p_0,\dots ,p_n)=(0,\dots
,0,p_n)$. Thus, we get $P=(p_0,\dots ,p_n,p_{n+1})=(0,\dots
,0,p_n,p_{n+1})\in L$. If $p_0=1$ then by assumption we know that
$P=(1,p_1,\dots ,p_n,p_{n+1}) \in Z(g_1,\dots ,g_{n-1},G)$. Since
$C_{\ell,m}$ is a s.t.c.i. on the hypersurfaces
$g_1=\cdots=g_{n-1}=0$ and $G=0$ it follows that $P=(1,p_1,\dots
,p_n,p_{n+1}) \in C_{\ell,m} \subset \overline{C}_{\ell,m}$.
\end{proof}

To get rid of $L$ in the intersection of the hypersurfaces
$f_1=\cdots=f_{n-1}=0$ and $F=0$, we modify the
$F={x_{1}}^{s_{1}}\cdots\,\,{x_{n}}^{s_{n}}-x_0^{\delta(m)-
\ell}x_{n+1}^{\ell}$ of the Lemma~\ref{projlemma}, as in the work
of Bresinsky (see \cite{bre}), for some special choice of
$f_1,\dots,f_{n-1}$. In this way we construct a new polynomial
$F^{*}$ from $F$ such that
$Z(f_1,\dots,f_{n-1},F^{*})=\overline{C}_{\ell,m}$, where $F^{*}$
is a polynomial of the form

\[ F^*=x_n^{\alpha}+x_0^{\beta}H(x_0,\dots ,x_{n+1}),\]
where $\beta$ is a positive integer.

Note that when  $x_0=0$, the vanishing of $F^{*}$ implies that
$x_n=0$. It follows from the last part of the proof of Lemma
\ref{projlemma} that this property of $F^{*}$ ensures that we have
a point at infinity, in the intersection of $f_1=\cdots=f_{n-1}=0$
and $F^{*}=0$, instead of a line.

The construction of $F^{*}$ can be described as follows. We first
assume that $f_i=x_i^{a_i}-x_0^{a_i-b_i}x_{n}^{b_i}=0$, where
$a_i>b_i$ are positive integers, for all $i=1,\dots ,n-1$. Let
$p=a_1\cdots a_{n-1}$ and $p_i=\frac{b_i}{a_i}p$, for $i=1,\dots
,n-1$. Take the $p$-th power of $F$ and for every occurrence of
$x_i^{a_i}$ substitute $x_0^{a_i-b_i}x_n^{b_i}$, for all
$i=1,\dots ,n-1$. Then we have
\begin{eqnarray*}
F^p&=&x_0^{\gamma}x_n^{\alpha}+ x_0^{\delta(m)-\ell}H(x_0,\dots
,x_{n+1}) \,\,\,\, mod(f_1,\dots,f_{n-1})\\
&=&x_0^{\gamma}[x_n^{\alpha}+x_0^{\delta(m)-\ell-\gamma}H(x_0,\dots
,x_{n+1})] \,\,\,\, mod(f_1,\dots ,f_{n-1})
\end{eqnarray*}
where $\gamma=\sum_{i=1}^{n-1} (p-p_i)s_{i}$,
$\alpha=ps_n+\sum_{i=1}^{n-1}p_i s_{i}$ and $H$ is a polynomial.
Letting
\[F^{*}(x_0,\dots ,x_{n+1})=x_n^{\alpha}+x_0^{\delta(m)-\ell-\gamma}H(x_0,\dots
,x_{n+1})\] we observe that
\begin{equation}\label{eqn}
F^p(x_0,\dots ,x_{n+1})=x_0^{\gamma}F^*(x_0,\dots ,x_{n+1})
\,\,\,\, mod(f_1,\dots ,f_{n-1}).
\end{equation}

Recall that $m$ is an element of the numerical semigroup generated
by $m_1,\dots,m_n$, i.e. $m=s_1m_1+\cdots+s_nm_n$ with
$s_1+\cdots+s_n=\delta(m)$. If $m$ is large enough that $s_n
> \ell+\sum_{i=1}^{n-1} (p-p_i-1)s_{i}$ (or equivalently
$\delta(m)-\ell-\gamma>0$) then $F^*$ is the required polynomial.
(Otherwise, $F^*$ may not be a polynomial.) Hence we conclude the
following

\begin{theorem}\label{projstci1} Let $p$, $p_i$, $f_i$ and $F^*$ be as above. Assume that $m$ is chosen so that $s_n > \ell+ \sum_{i=1}^{n-1} (p-p_i-1)s_{i}$.
Then, for all $\ell<\delta(m)$ with $gcd(\ell,m)=1$, the nice
extensions $\overline{C}_{\ell,m} \subset \mathbb{P}^{n+1}$ are
s.t.c.i. on $f_1=\cdots=f_{n-1}=0$ and $F^*=0$.

\end{theorem}

\begin{proof} We will show that $\overline{C}_{\ell,m}$ is a
s.t.c.i. on $f_1=\cdots=f_{n-1}=0$ and $F^*=0$. To do this, take a
point $P=(p_0,\dots ,p_{n+1}) \in \overline{C}_{\ell,m}$. Then,
$F(P)=0$ and $f_i(P)=0$, for all $i=1,\dots ,n-1$, since
$Z(f_1,\dots ,f_{n-1},F)=\overline{C}_{\ell,m} \bigcup L$, by
Lemma \ref{projlemma}. From equation (\ref{eqn}) it follows that
$F^*(P)=0$ or $p_0=0$. Since $P$ is a point on the monomial curve
$\overline{C}_{\ell,m}$, it can be parameterized as follows:
\[(u^{m},u^{m-\ell m_{1}}v^{\ell m_{1}},\dots
,u^{m-\ell m_{n}}v^{\ell m_{n}},v^{m})
\]

So if $p_0=0$, we get $u=0$ and thus $p_i=0$, for all $i=1,\dots
,n$. Therefore $P=(0,\dots ,0,1)$ and hence $F^*(P)=0$ in any
case.

Conversely, let $P=(p_0,\dots ,p_{n+1}) \in Z(f_1,\dots
,f_{n-1},F^*)$. If $p_0=0$, then $p_i=0$ by $f_i(P)=0$, for all
$i=1,\dots ,n-1$. Since $\delta(m)-\ell-\gamma>0$, we have $p_n=0$
by $F^*(P)=0$. Thus $P=(0,\dots ,0,1)$ which is always on the
curve $\overline{C}_{\ell,m}$. If $p_0=1$ then $C$ is a s.t.c.i.
on the hypersurfaces given by $g_i=x_i^{a_i}-x_{i+1}^{b_i}=0$, for
$i=1,\dots ,n-1$, by the assumption. Hence, Theorem \ref{afinstci}
implies that $C_{\ell,m}$ is a s.t.c.i. on $g_1=\cdots=g_{n-1}=0$
and $G={x_{1}}^{s_{1}}\cdots\,\,{x_{n}}^{s_{n}}-x_{n+1}^{\ell}=0$.
Thus $P=(1,p_1,\dots,p_{n+1}) \in C_{\ell,m} \subset
\overline{C}_{\ell,m}$.
\end{proof}

\begin{remark} The \emph{nice} extensions in Theorem \ref{projstci1} can also be shown to be s.t.c.i. by using
\cite[Theorem~3.4]{thoma2}. But to show that the hypotheses of
\cite[Theorem~3.4]{thoma2} are satisfied by these extensions is
much more difficult than the proof here. As a byproduct we also
constructed here the hypersurface $F^*=0$ on which these
\emph{nice} extensions are s.t.c.i.
\end{remark}

\begin{example} We start with $\overline{C}=\overline{C}(3,4,6)\subset
\mathbb{P}^{3}$. Let $\ell=1$ and $m=6s+7$, for some positive
integer $s$. Then $\delta(m)=s+2$, $s_1=s_2=1$ and $s_3=s$. Thus
we get the nice extensions
$\overline{C}_{1,6s+7}=\overline{C}(3,4,6,6s+7)\subset
\mathbb{P}^{4}$. Since $\Delta_1=gcd(4,6,6s+7)=1$,
$\Delta_2=gcd(3,6,6s+7)=1$ and $\Delta_3=gcd(3,4,6s+7)=1$ it
follows from Corollary \ref{gluingcor1} that these curves can not
be obtained by gluing. Using the software Macaulay, it is easy to
see that the ideal of $\overline{C}_{1,6s+7}$ is minimally
generated by the polynomials
\begin{eqnarray*} f_1&=&x_1^2-x_0x_3,\\ f_2&=&x_2^3-x_0x_3^2, \\
f_3&=&x_3^{s+3}-x_0^{s-1}x_1x_2^2x_4 \\
f_4&=&x_2x_3^{s+1}-x_0^{s}x_1x_4, \\
f_5&=&x_1x_3^{s+2}-x_0^{s}x_2^2x_4 \\
F&=&x_1x_2x_3^{s}-x_0^{s+1}x_4.
\end{eqnarray*}
Since $\overline{C}(3,4,6)\subset \mathbb{P}^{3}$ is a s.t.c.i. on
the surfaces $f_1=0$ and $f_2=0$, it follows from Theorem
\ref{projstci1} that $\overline{C}_{1,6s+7}$ is a s.t.c.i. on
$f_1=0$, $f_2=0$ and
$$F^*=x_3^{6s+7}-6x_0^{s-1}x_1x_2^2x_3^{5s+4}x_4+
15x_0^{2s}x_2x_3^{4s+4}x_4^2-20x_0^{3s}x_1x_3^{3s+3}x_4^3+$$
$$+15x_0^{4s}x_2^2x_3^{2s+1}x_4^4-6x_0^{5s}x_1x_2x_3^sx_4^5+
x_0^{6s+1}x_4^6=0$$ provided that $s>2$.

\end{example}

\section{Variations and consequences of the main results}

In this section, we give some  consequences of Theorem
\ref{projstci} and hence all the notation is as in that theorem.
We also include some theorems about \emph{nice} extensions of
projective monomial curves that are variations of Theorem
\ref{projstci1}.

\subsection{Consequences of Theorem \ref{projstci}}

Since arithmetically Cohen-Macaulay monomial curves are s.t.c.i.
in $ \mathbb{P}^3$ (see \cite{rv-ACM}), we get the following
corollary as a consequence of Theorem \ref{projstci}.

\begin{corollary}\label{corprojmain} Let $\overline{C}(m_1,m_2,m_3)$ be an arithmetically
Cohen-Macaulay monomial curve in $\mathbb{P}^3$. Let $m=s_3m_3$,
$gcd(\ell,m)=1$ and $\delta(m)=s_3>\ell$. Then the nice extensions
$\overline{C}_{\ell,s_3m_3}=\overline{C}(\ell m_1,\ell m_2,\ell
m_3,s_3m_3)$ are all s.t.c.i. in $\mathbb{P}^4$. \mbox{} \hfill
$\Box$
\end{corollary}

\begin{remark} There are very few examples of s.t.c.i. monomial curves in $\mathbb{P}^n$, where $n>3$. We
know that rational normal curve $\overline{C}(1,2,\dots ,n)$ is a
s.t.c.i. in $\mathbb{P}^{n}$, for any $n>0$, (see
\cite{rv-normal,thoma2}). Applying Theorem \ref{projstci} to
$\overline{C}(1,2,\dots ,n) \subset \mathbb{P}^{n}$, we can
produce infinitely many new examples of s.t.c.i. monomial curves
in $\mathbb{P}^{n+1}$:
\end{remark}

\begin{corollary}
For all positive integers $\ell$, $n$ and $s$ with
$gcd(\ell,sn)=1$, the monomial curves
$\overline{C}(\ell,2\ell,\dots ,n\ell,sn) \subset
\mathbb{P}^{n+1}$ are s.t.c.i.
\end{corollary}

\begin{proof} Let $m=sn$. Clearly $\delta(m)=s$. If $s \leq \ell$,
then the curves
$\overline{C}_{\ell,m}=\overline{C}(\ell,2\ell,\dots ,n\ell,sn)
\subset \mathbb{P}^{n+1}$ are bad extensions of
$\overline{C}(1,2,\dots,n)\subset \mathbb{P}^n$. Hence they are
s.t.c.i. by Corollary \ref{badextensions}. If $s>\ell$, then these
curves are nice extensions of $\overline{C}(1,2,\dots,n)\subset
\mathbb{P}^n$. Therefore they are s.t.c.i. by Theorem
\ref{projstci}.
\end{proof}

In \cite{mt}, all complete intersection (i.t.c.i.) lattice ideals
are characterized by gluing semigroups. But, for a given
projective monomial curve it is not easy to find two subsemigroups
whose ideals are complete intersection. So, as another application
of Theorem \ref{projstci} we can produce infinitely many i.t.c.i.
monomial curves:

\begin{proposition}\label{projitci} If \,$\overline{C}\subset\mathbb{P}^{n}$ is an i.t.c.i., then the nice extensions $\overline{C}_{\ell,s_nm_n}\subset\mathbb{P}^{n+1}$
are i.t.c.i. for all positive integers $\ell$ and $s_n$ with
$s_n>\ell$, $gcd(\ell,s_nm_n)=1$.
\end{proposition}

\begin{proof} Since $\overline{C}$ is a s.t.c.i. on the binomial hypersurfaces
$f_1=\cdots=f_{n-1}=0$, it follows from Theorem \ref{projstci}
that $\overline{C}_{\ell,s_nm_n}$ is a s.t.c.i. on
$f_1=\cdots=f_{n-1}=0$ and $F(x_0,\dots
,x_{n+1})=x_{n}^{s_n}-x_{0}^{s_n-\ell}x_{n+1}^{\ell}=0$. Since
these are all binomial, the monomial curves
$\overline{C}_{\ell,s_nm_n}$ are i.t.c.i. on the same
hypersurfaces, by \cite[Theorem~4]{bmt}.

\end{proof}

\begin{corollary}\label{itci} The monomial curves $\overline{C}(\ell m_1,\ell m_2,s_2m_2)$ are i.t.c.i. in
$\mathbb{P}^{3}$, for all positive integers $m_1,m_2,\ell$ and
$s_2$ with $s_2>\ell$, $gcd(\ell,s_2m_2)=1$.
\end{corollary}

\begin{proof} Let $m=s_2m_2$. Then $\delta(m)=s_2$ and $\overline{C}_{\ell,m}=\overline{C}(\ell m_1,\ell m_2,s_2m_2)$ is a nice extension of
$\overline{C}(m_1,m_2)$, by the assumption $s_2>\ell$. Since
$\overline{C}(m_1,m_2)$ is an i.t.c.i. on
$x_1^{m_2}-x_0^{m_2-m_1}x_2^{m_1}=0$, it follows from Proposition
\ref{projitci} that the nice extensions $\overline{C}(\ell
m_1,\ell m_2,s_2m_2)$ are i.t.c.i. on
$x_1^{m_2}-x_0^{m_2-m_1}x_2^{m_1}=0$ and
$x_2^{s_2}-x_0^{s_2-\ell}x_3^{\ell}=0.$
\end{proof}
To produce infinitely many examples of i.t.c.i. curves, our method
starts from just one i.t.c.i. curve, whereas semigroup gluing
method produces only one example starting from one i.t.c.i.. The
following example illustrates this point.
\begin{example}\label{exampleitci} From Corollary \ref{itci}, we know that $\overline{C}(1,2,4)$ is an i.t.c.i. on
\[ f_1=x_1^2-x_0x_2=0\; \text{~and~} \; f_2=x_2^2-x_0x_3=0. \]
Take two positive integers $\ell$ and $s$ with $s>\ell$,
$gcd(\ell,4s)=1$. Then the monomial curves
$\overline{C}(\ell,2\ell,4\ell,4s)\subset \mathbb{P}^4$ are nice
extensions of $\overline{C}(1,2,4)\subset \mathbb{P}^3$. Thus, by
Proposition \ref{projitci}, the monomial curves
$\overline{C}(\ell,2\ell,4\ell,4s)$ are i.t.c.i. on
\[ f_1=x_1^2-x_0x_2=0, \; f_2=x_2^2-x_0x_3=0 \; \text{~and~} \;
F=x_3^s-x_0^{s-\ell}x_4^{\ell}=0.
\] The nice extensions
$\overline{C}(\ell,2\ell,4\ell,4s)$ can also be obtained by gluing
subsemigroups generated by $T_1=\{(4s-\ell,\ell)\}$ and
$T_2=\{(4s,0),(4s-2\ell,2\ell),(4s-4\ell,4\ell),(0,4s)\}$. But, in
this case one has to know that $\overline{C}(\ell,2\ell,2s)$ is an
i.t.c.i. for each $\ell$ and $s$. In other words, starting with
the fact that $\overline{C}(1,2,4)$ is an i.t.c.i., gluing method
can only produce $\overline{C}(1,2,4,8)$ as an i.t.c.i. monomial
curve.
\end{example}

\subsection{Variations of Theorem \ref{projstci1}} Recall that our method starts with a monomial curve
$\overline{C}=Z(f_1,\dots,f_{n-1})$ in $\mathbb{P}^n$ and produces
infinitely many nice extensions
$\overline{C}_{\ell,m}=Z(f_1,\dots,f_{n-1},F^*)$ in
$\mathbb{P}^{n+1}$. Since the construction of $F^*$ depends on the
choice of $f_1,\dots, f_{n-1}$, it is possible to start with
another curve $\overline{C}=Z(f_1,\dots,f_{n-1})$ in
$\mathbb{P}^n$ and obtain new families of nice extensions. In this
section we provide two examples of this sort. For instance, if we
assume that $\overline{C}$ is a s.t.c.i. on the hypersurfaces
$f_i=x_i^{a_i}-x_0^{a_i-b_i}x_{i+1}^{b_i}=0$, where $a_i > b_i$
are positive integers, $i=1,\dots ,n-1$, then under some suitable
conditions we obtain other families of s.t.c.i. nice extensions.
Let $p=a_1\cdots a_{n-1}$, $q_0=b_1 \cdots b_{n-1}$ and
$q_i=a_1\cdots a_ib_{i+1}\cdots b_{n-1}$, $i=1,\dots ,n-2$. The
first variation is the following
\begin{theorem}\label{projstci3} Let $p,q_0,\dots,q_{n-2}$ be as above. For all $m$ which give rise to $s_n
>\ell+\sum_{i=0}^{n-2}(p-q_i-1)s_{i+1}$ and for all $\ell$ with $\ell<\delta(m)$ and $gcd(\ell,m)=1$, the nice
extensions $\overline{C}_{\ell,m}\subset\mathbb{P}^{n+1}$ are
s.t.c.i. on $f_1=\cdots=f_{n-1}=F^*=0$. \hfill $\Box$
\end{theorem}

Now, we give another variation where $m=s_im_i+s_jm_j$, for
$i,j\in \{1,\dots,n\}$. For the notational convenience we take
$i=1$ and $j=n$.

\begin{theorem}\label{projstci2} Let $\overline{C}\subset \mathbb{P}^{n}$ be a s.t.c.i.
on the hypersurfaces given by
\begin{eqnarray*}f_1&=&x_1^{a}-x_0^{a-b}x_{n}^{b}=0\\
f_i&=&x_i^{a_i}+x_0^{b_i}A(x_1,\dots ,x_n)+x_{1}^{c_i}B(x_2,\dots
,x_n)=0,
\end{eqnarray*} where $a$, $b$, $a-b$, $a_i$, $b_i$, and $c_i$ are positive
integers, for $i=2,\dots ,n-1$, $A$ and $B$ are some polynomials.
For all $m$ which give rise to $s_n > \ell+(a-b-1)s_1$ and for all
$\ell$ with $\ell<\delta(m)$ and $gcd(\ell,m)=1$, the nice
extensions $\overline{C}_{\ell,m}\subset\mathbb{P}^{n+1}$ are
s.t.c.i. on $f_1=\cdots=f_{n-1}=F^*=0$. \hfill $\Box$
\end{theorem}

\section*{Acknowledgements} The author would like to extend his sincere thanks to F. Arslan,
\"{O}. Ki\c{s}isel, M. Morales, S. Sert\"{o}z and A. Thoma for
their numerous comments, and to the anonymous referee for his
suggestions which improved the final presentation of the paper.


\bigskip

\begin{thebibliography}{99}

\bibitem{pf}F. Arslan, P. Mete, \textit{Hilbert functions of Gorenstein monomial curves},
Proc. Amer. Math. Soc. \textbf{135} (2007) 1993-2002.

\bibitem{bmt}M. Barile, M. Morales, A.Thoma, \textit{Set-theoretic complete intersections on binomials},
Proc. Amer. Math. Soc. \textbf{130} (2002) 1893-1903.

\bibitem{bre}H. Bresinsky, \textit{Monomial space curves in $\mathbb{A}^3$ as
set-theoretic complete intersection}, Proc. Amer. Math. Soc.
\textbf{75} (1979) 23-24.

\bibitem{eis}D. Eisenbud, E.G. Evans, \textit{Every Algebraic Set in n-space is the intersection of n Hypersurfaces},
Inventiones Math. \textbf{19} (1973) 107-112.

\bibitem{eto2}K. Eto, \textit{Set-theoretic complete intersection lattice ideals in
monoid rings}, Journal of Algebra \textbf{299} (2006) 689-706.

\bibitem{kat} A. Katsabekis, \textit{Projection of cones and the arithmetical rank of toric varieties},
Journal of Pure and Applied Algebra, \textbf{199} (2005) 133-147.

\bibitem{mac}D. Bayer and M. Stillman, \textit{Macaulay, A system for computations in
algebraic geometry and commutative algebra}, 1992.

\bibitem{moh}T.T. Moh, \textit{Set-theoretic complete intersections},
Proc. Amer. Math. Soc. \textbf{94} (1985) 217-220.

\bibitem{mor}M. Morales, \textit{Noetherian Symbolic Blow-Ups}, Journal
of Algebra \textbf{140} (1991) 12-25.

\bibitem{mt}M. Morales and A. Thoma, \textit{Complete intersection lattice ideals}, Journal
of Algebra \textbf{284} (2005) 755-770.

\bibitem{rv-normal}L. Robbiano, G. Valla, \textit{On set-theoretic complete
intersections in the projective space}, Rend. Sem. Mat. Fis.
Milano LIII (1983) 333-346.

\bibitem{rv-ACM}L. Robbiano, G. Valla, \textit{Some curves in $\mathbb{P}^{3}$ are
set-theoretic complete intersections}, in: Algebraic Geometry-Open
problems, Proceedings Ravello 1982, Lecture Notes in Mathematics,
Vol \textbf{997} (Springer, New York, 1983) 391-346.

\bibitem{ros}J.C. Rosales, \textit{On presentations of subsemigroups of $\mathbb{N}^n$},
Semigroup Forum \textbf{55} (1997) 152-159.

\bibitem{thoma2}A. Thoma, \textit{On the set-theoretic complete intersection
problem for monomial curves in $\mathbb{A}^{n}$ and
$\mathbb{P}^{n}$}, Journal of Pure and Applied Algebra,
\textbf{104} (1995) 333-344.

\bibitem{thoma3}A.Thoma, \textit{Affine semigroup rings and monomial varieties},
Communications in Algebra \textbf{24(7)} (1996) 2463-2471.

\bibitem{thoma4}A.Thoma, \textit{Construction of set-theoretic complete intersections via semigroup gluing},
Contributions to Algebra and Geometry \textbf{41(1)} (2000)
195-198.


\end{thebibliography}
\end{document}